\documentclass[final]{siamltex}

\usepackage{amssymb}
\usepackage{amsmath}
\usepackage{graphicx}
\usepackage{subfigure}
\usepackage{color}

   \newcommand{\RR}{{\mathbb R}}

\newcommand{\veps}{\varepsilon}

\begin{document}
 
\title{Implicit-Explicit Runge-Kutta schemes for hyperbolic systems
        with stiff relaxation\\ and applications}

\author{Sebastiano Boscarino\thanks{Department of
Mathematics and Computer Science, University of Catania, Via
A.Doria 6, 95125 Catania, Italy. ({\tt boscarino@dmi.unict.it, russo@dmi.unict.it})}\and
%
Giovanni Russo$^\star$}

\maketitle

\begin{abstract} In this paper we give an overview of Implicit-Explicit Runge-Kutta
schemes applied to 
hyperbolic systems with stiff relaxation. In particular, we focus on
some recent results on the uniform accuracy for hyperbolic systems
with stiff relaxation \cite{BR2009}, and hyperbolic system with
diffusive relaxation \cite{BR1,BPR,BFR}. In the latter case, we
present an original application to a model problem arising in Extended Thermodynamics.
\end{abstract}



\pagestyle{myheadings}
\thispagestyle{plain}
\markboth{S. BOSCARINO AND G. RUSSO}{} 


\section{Introduction}

Many physical models are described by hyperbolic systems with
relaxation of the form
\begin{equation}
   \partial_t U + \partial_x F(U) = \frac1{\varepsilon} R(U),\quad
     x \in \RR,
    \label{eq:hyp-rel}
\end{equation}
with $U=U(x,t) \in \RR^N$, $F:\RR^N\to\RR^N$. Such systems are said hyperbolic if the Jacobian matrix $F'(U)$ has real eigenvalues and
a basis of eigenvectors $\forall\,U\in \RR^N$. Usually, the parameter $\varepsilon$ is called the relaxation
time, which is small in many physical situations.
Here we use the term relaxation in the sense of Whitham \cite{Whith}
and Liu (\cite{Liu2}), which in practice means that if $\veps\to
0$, the system formally relaxes to a quasilinear hyperbolic system
with a smaller number of dimensions. Chen, Levermore, and Liu \cite{Liu1} provide
the proper condition that ensures that the solution of the relaxation
system actually converges to the solution of the {\em relaxed\/}
system. 

Typical examples of such systems are: gas dynamics with chemical reactions,  
shallow water with friction, 
discrete kinetic models, 
extended thermodynamics, 
hydrodynamical models for semiconductors,  
traffic flow models, 
granular gases (see \cite{typeA} and references therein). 

 A simple prototype example of relaxation system is given by
\begin{eqnarray*}
\partial_t u + \partial_x v &=&  0,\\
\partial_t v + \partial_x p(u) &=& -\frac{1}{\varepsilon} (v-f(u)),
\end{eqnarray*}
which corresponds to $U=(u,v)$,
$F(U)=(v,p(u))$,  $R(U)=(0,f(u)-v)$. 
As $\varepsilon \to 0$ we get, formally, the local equilibrium $v=f(u)$
while $u$ satisfies the conservation equation
\[
  \partial_t u + \partial_x f(u)=0.
\]
In \cite{Liu1} the authors proved
that the solution $u$ actually converges to the solution of the relaxed equation
if the characteristic speed of the relaxed equation is contained in the
interval identified by the speed of the original system,i.e. $p'(u) \ge (f'(u))^2$, i.e.
the \emph{subcharacteristic condition}.

The most commonly used approach for the numerical solution of hyperbolic
system with relaxation is based on the Method Of Line (MOL).
First we discretize the system in space, leading to a large system of ODEs defined
on a grid. The semi discrete scheme should be high resolution shock
capturing, which provide correct shock location without numerical oscillations.
Among space discretization techniques we mentioned several
possibilities: Finite Volume (FV), Finite Difference (FD),
Discontinuous Galerkin (DG).
Method of lines based on conservative finite difference is the
simplest choice for the construction of high order schemes in space
and time \cite{Shu, typeA}.
For example, in one space dimension, the scheme reads:
\[
  \frac{du_j}{dt} = -\frac{\hat{f}_{j + 1/2}-\hat{f}_{j - 1/2}}{dx} - g(u_j)
\]
with 
\[
   \hat{f}_{j + 1/2} = \hat{f}_{j + 1/2}^+(x_{j + 1/2}^-) + \hat{f}_{j+ 1/2}^-(x_{j+1/2}^+).
\]
The numerical flux $\{\hat{f}_{j + 1/2}^\pm(x)\}$ being reconstructed from the fluxes
$f^\pm(x_j)$, which in turn split the analytical flux: $f = f^+ + f^-$, $\lambda(\nabla f^+) \ge 0$, $\lambda(\nabla f^+) \le 0$. High
order reconstruction can be obtained, for example, by ENO or WENO
reconstruction from cell averages to pointwise values,
\[
    \{ f^\pm_j\}   \xrightarrow[\phantom{WENO}]{WENO}  \hat{f}_j^\pm(x_{j\pm1/2})
\]
Since source term $g(u_j )$ is computed \emph{pointwise} then the various cells are not
coupled at the level of the source, and the implicit equations in each cell are
independent from each other.

Applying MOL to hyperbolic system with relaxation, the PDEs become
 a system of ODEs of the form
\begin{equation}
  u'  =  f(u) + \frac1{\varepsilon}g(u),
  \label{ODE}
\end{equation}
with initial vector $u_0 = (U(x_1,t_0),\cdots, U(x_N, t_0))^T$, where
$\{x_i\}_{i = 1}^{N}$ denote the spatial computational
mesh. The solution at time $t$ is $u(t) = (u(t_1),
u_2(t),\cdots,u_N(t))^T$ where $u_i(t) \approx U(x_i,t)$. The term
$f(u)$ represents the discretization of the convective derivative
term, $- \partial_x F(U)$, while $g(u)$ represents the discrete
approximation of the source term, $G(U)$, on the grid nodes (and
possibly the boundary conditions). Then a suitable time
integrator is used to solve ODEs.

In most cases $f(u)$ is non stiff and non linear while
$\frac{1}{\varepsilon}g(u)$ contains the stiffness, so we look for 
numerical schemes which are explicit in $f$ and implicit in
$g$.  In particular it is essential that the numerical scheme is
accurate for $ \varepsilon \to 0$ (possibly also for intermediate
regimes of such parameters, i.e.\ when $\varepsilon$ is not too small).
Moreover some stability restrictions are required, i.e. for the
convection term $\Delta t \leq \rho(\nabla_u F)\Delta x$ (CFL
condition). The stiff term has to be treated implicitly to avoid
restrictions $\Delta t \leq C \varepsilon$.

IMEX Runge-Kutta methods represents a very effective tool to guarantee
the simplicity of the explicit treatment of the non-stiff term $f(u)$
and to avoid time restriction because of the stiffness in the source
term $g(u)$.

An Implicit-Explicit (IMEX) Runge-Kutta scheme applied to system
(\ref{ODE}) takes the form 
\begin{eqnarray*}
 Y_i & = & y_0 + h \sum_{j=1}^{i-1} \tilde{a}_{ij} f(t_0+\tilde{c}_jh, Y_j) +
                 h \sum_{j=1}^i a_{ij} \frac{1}{\varepsilon} g(t_0+c_jh, Y_j),
                                                                 \label{eq:RKEI1}\\
 y_1 & = & y_0 + h \sum_{i=1}^{s} \tilde{b}_{i} f(t_0+\tilde{c}_ih, Y_i)  +
                 h \sum_{i=1}^s b_{i} \frac{1}{\varepsilon} g(t_0+c_ih, Y_i).
                                                                 \label{eq:RKEI2}
\end{eqnarray*}
where $\tilde{A} = (\tilde{a}_{ij})$, $\tilde{a}_{ij}=0$, $j \geq i$
and $A = (a_{ij})$ are $s \times s$ (lower triangular) matrices and
$\tilde{c}, \tilde{b}, c, b \in \RR^s$, coefficient vectors. A
classical representation of a IMEX R-K method is given by
\[
\hbox{Double Butcher {\em tableau\/}:} \quad
\begin{array}{c|c}
              \tilde{c} & \tilde{A} \\
              \hline \\
              & \tilde{b}^T
\end{array} \qquad
\begin{array}{c|c}
              c & A \\
              \hline \\
              & b^T
\end{array}.
\]
We restrict to consider IMEX schemes in which the implicit part is a
diagonally implicit Runge-Kutta (DIRK). Besides it simplicity, this
will ensure that $f$ is always evaluated explicitly.

We can classify each IMEX Runge-Kutta scheme by considering the different
structures of the matrix $A = (a_{ij})_{i,j=1}^s$, of the
implicit scheme:
\begin{itemize}
\item  (Methods of Type A) The matrix A is invertible.
\item (Methods of Type CK)
\[
      A = \left(\begin{array}{ll} 0 & 0\\
                        a  & \hat{A}
          \end{array}\right)
\]         
The submatrix $\hat{A}$ is invertible. 
\end{itemize}
CK methods with $a=0$ are called ARS methods \cite{ARS}. 
Type A methods are somehow more difficult to construct, but easier to
analyze than methods of type CK \cite{CK} or ARS. 

The rest of the paper is organized as follows.
Section \ref{Sec1} review some recent results on the development of high-order
implicit-explicit (IMEX) Runge-Kutta (R-K) schemes suitable for time-dependent
partial differential systems \cite{BR2009}. In section \ref{Sec2} we discuss hyperbolic 
systems with stiff diffusive relaxation. The last section is devoted
to some applications to some models of diffusive relaxation, which
confirm practice the advantageous effects of the approaches introduced
the earlier sections.  
In particular, Sec.~\ref{R13} is devoted to an original application of
IMEX-I schemes without parabolic restriction to a one dimensional model problem
arising in the context of Extended Thermodynamics.

\section{On the uniform accuracy of IMEX RungeÐKutta
schemes and applications to hyperbolic systems
with relaxation.}\label{Sec1}

Usually, under-resolved numerical schemes may yield spurious numerical
solutions that are unphysical. Other times, in the case of hyperbolic
systems with stiff terms, high order schemes may reduce to lower order
when the time step fails to resolve the small relaxation time.

IMplicit-EXplicit (IMEX) Runge-Kutta (R-K) schemes have been widely
used for the time evolution of hyperbolic partial differential
equations but some of the schemes existing in literature do not exhibit
uniform accuracy with respect to the relaxation time. 
Classical
high-order IMEX R-K schemes fail to maintain the high-order accuracy
in time in the whole range of the relaxation time and in particular in
the asymptotic limit $\veps\to 0$.

In \cite{BR2009} we developed new IMEX R-K schemes for hyperbolic
systems with relaxation that present better uniform accuracy than the
ones existing in the literature and in particular produce good
behavior with high order accuracy in the asymptotic limit, i.e.\ when
$\varepsilon$ is very small. In particular, these schemes are able to
handle the stiffness of the system (\ref{eq:hyp-rel}), in a whole
range of the relaxation time.

The schemes are obtained by imposing new additional conditions on
their coefficients, in order to guarantee better accuracy over a wide
range of the relaxation time.  Following the same technique proposed in \cite{HWVol2}, the additional
conditions are obtained by performing an asymptotic expansion of the
exact and numerical solution in the small parameter $\veps$
(Hilbert expansion), and by matching the two solutions to various
order in $\veps$, \cite{Boscarino2009}. 

The construction of a high-order accurate
IMEX R-K scheme is obtained by imposing the extra order conditions, that ensure 
the agreement between exact and numerical solution up to  a given order in $\varepsilon$.
The scheme, called BHR(5,5,3),  
is presented in \cite{Boscarino2009, BR2009}. 
Numerical tests on several ordinary differential
systems and hyperbolic systems with relaxation term present 
better behavior for the new scheme BHR(5,5,3)
over other IMEX R-K methods previously existing in literature \cite{ARS, CK, typeA}.
For example, by imposing the additional order conditions to the zeroth-order in $\varepsilon$, the classical ARS(4,3,4) scheme
can be modified (hereafter called Mod-ARS(3,4,3)), imposing its accuracy in the algebraic variable.
 Furthermore, by imposing conditions to  terms up to fist order in $\varepsilon$
and we constructed scheme RHR(5,5,3), a third order five stage scheme.

\begin{figure}
\begin{center}
  \includegraphics[width=0.49\textwidth]{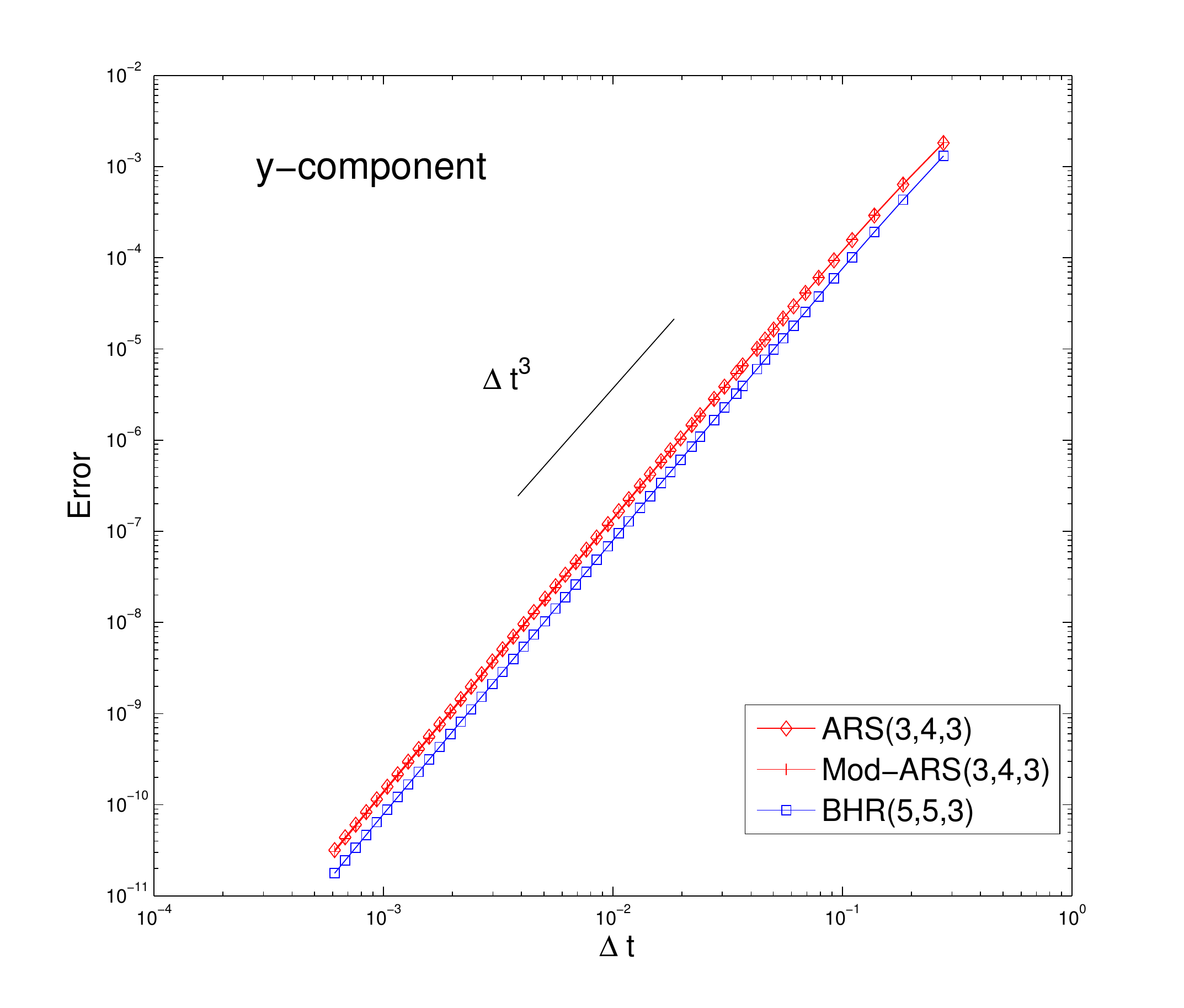}
  \includegraphics[width=0.49\textwidth]{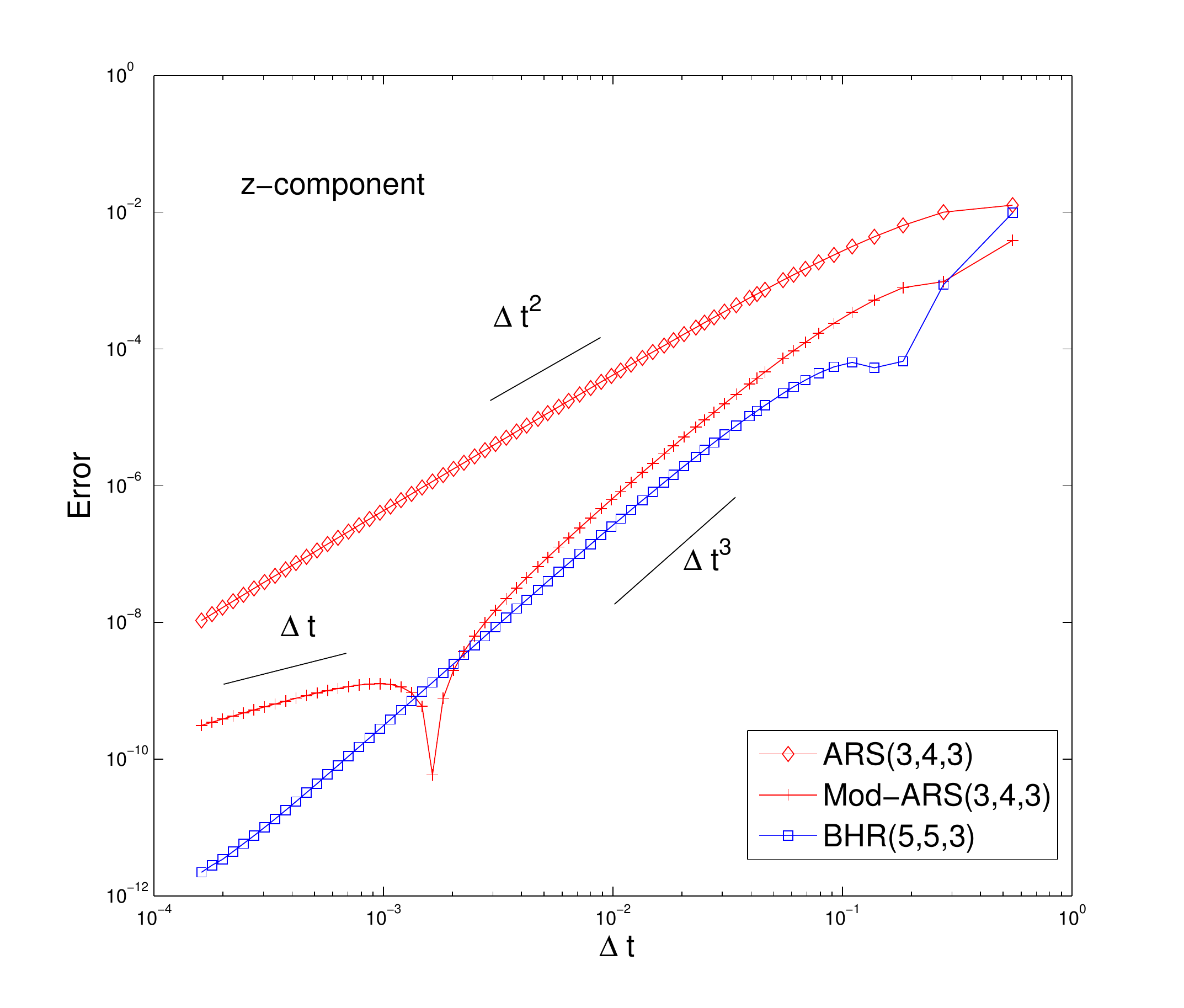}
  \caption{Global error versus the stepsize in the Van der Pol equation calculated with $\varepsilon = 10^{-6}$.}\label{fig1}
  \end{center}
\end{figure}

The construction of this type of IMEX R-K scheme is motivated by the
order reduction of classical IMEX schemes observed when applying them
to several stiff systems. An example of such behavior is illustrated
in Figure \ref{fig1}, where the classical Van rer Pol equation is solved by ARS(3,4,3), Mod-ARS(3,4,3) and BHR(5,5,3) schemes 
derived in \cite{ARS, CK, typeA, Boscarino2009, BR2009},
\begin{eqnarray}
\begin{array}{l}
y' = z,\\
\varepsilon z' = (1 - y^2)z-y,
\end{array}
\end{eqnarray}
 (for details of this problem
and its initial conditions see, for example, \cite{HWVol2}).
The global error behaves like $C \Delta t^r$ with $r$ the slope of the
straight line and $C$ is a constant. 
We observe that, while classical schemes, as ARS(3,4,3), are able to maintain the classical
order of accuracy in the differential variable $y$, they lose accuracy
in the algebraic variable $z$. 
BHR(5,5,3) method exhibits the better error estimate with respect to ARS(3,4,3) and Mod-ARS(3,4,3) schemes and no order reduction appears when $\varepsilon$ is very small.

Concerning hyperbolic systems with stiff relaxation  we report here a
numerical test the Broadwell model equations

\begin{eqnarray}
\begin{array}{l}
\rho_t + m_z = 0,\\
m_t + z_x = 0,\\
z_t + m_x = \frac{1}{\varepsilon}(\rho^2 + m^2 -2\rho z)
\end{array}
\end{eqnarray}

(for details see
\cite{BR2009}), which, in one space dimension, is a $3\times 3$ semilinear hyperbolic system that, in in
the relaxed limit, becomes a quasilinear hyperbolic system for the two
two differential variables ($\rho$ and $m$), while $z$ becomes a
function of the other two variables. 

\begin{figure}
\begin{center}
  \includegraphics[width=0.49\textwidth]{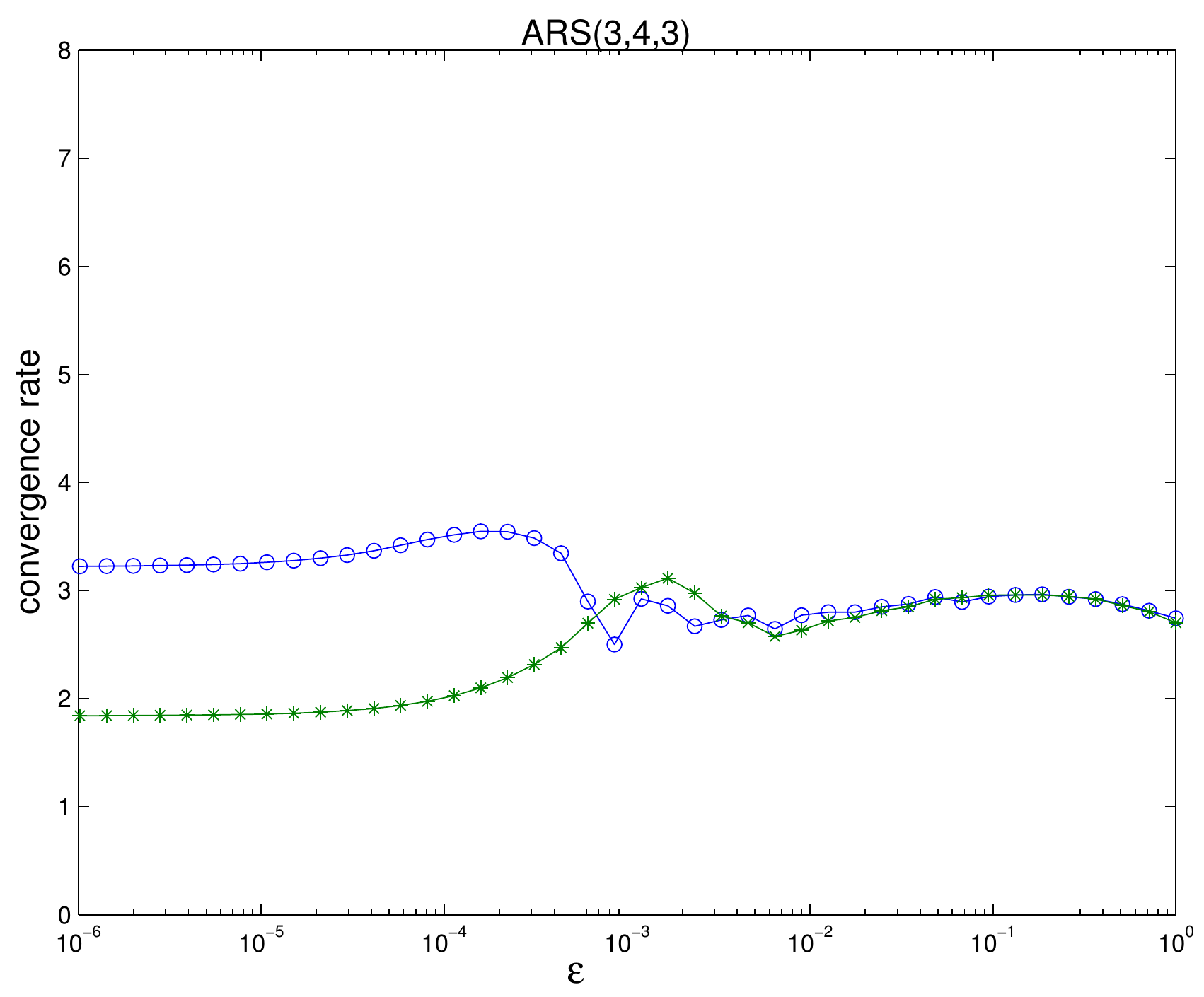}
   \includegraphics[width=0.49\textwidth]{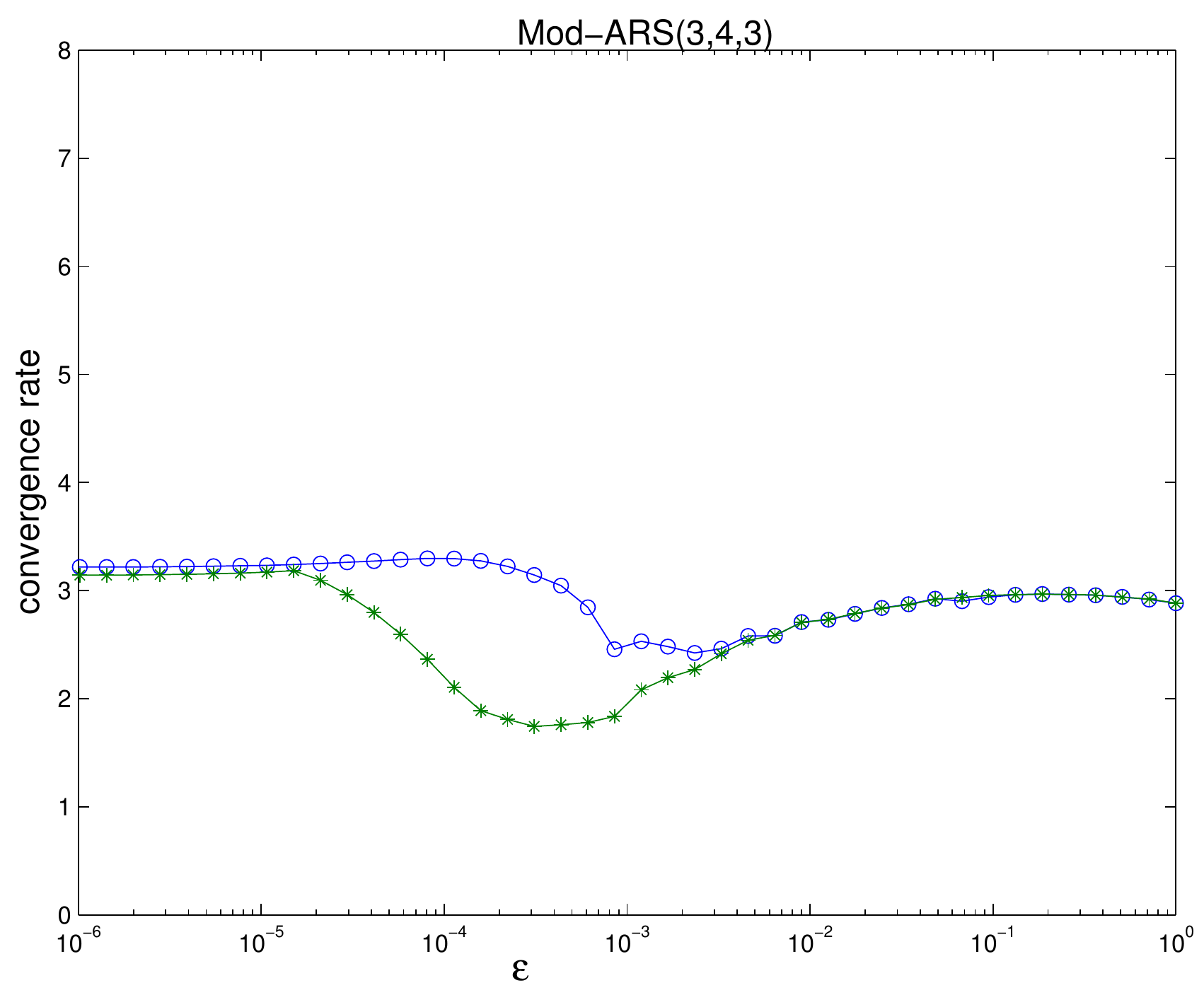}
   \includegraphics[width=0.49\textwidth]{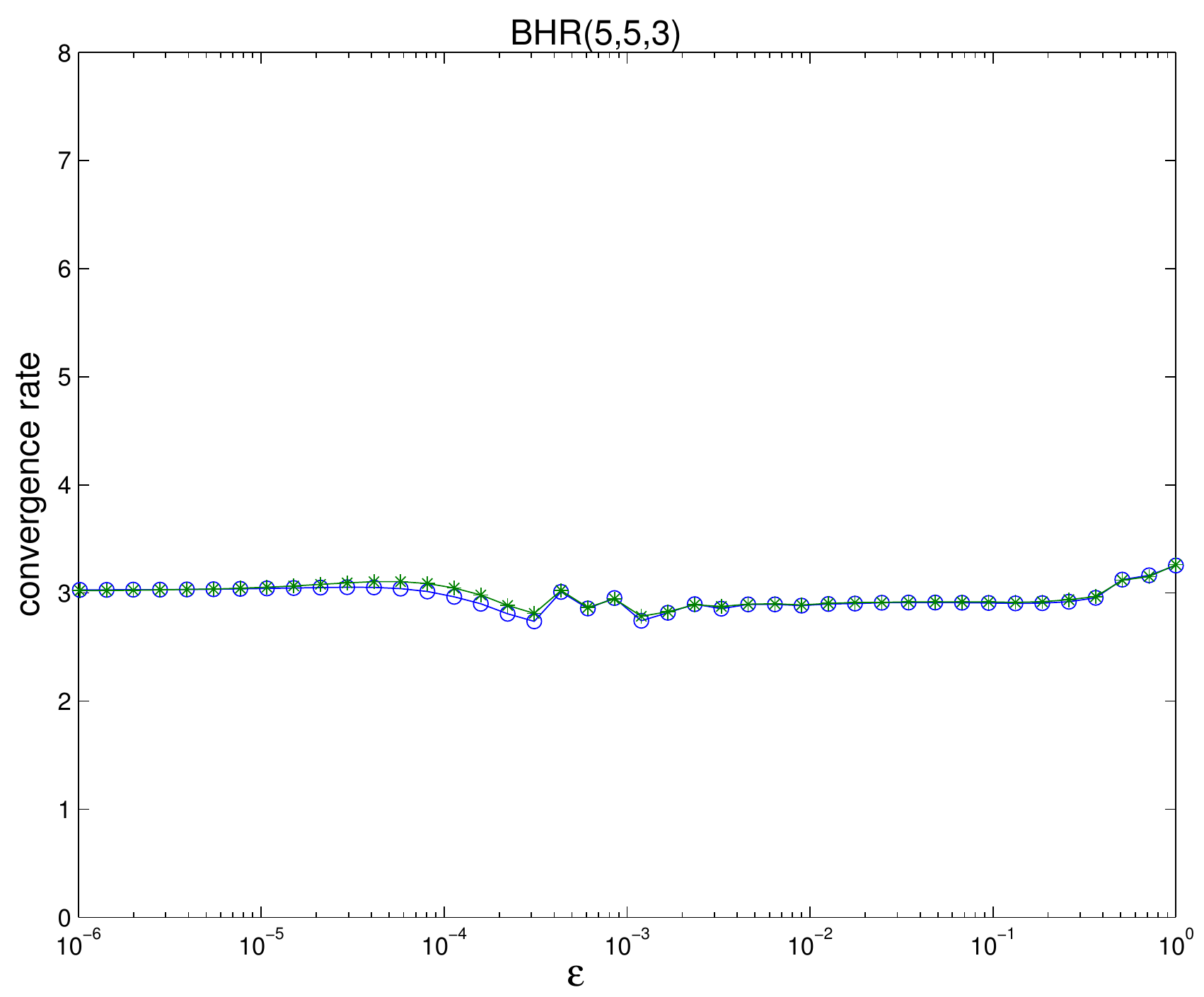}
  \caption{Convergence rate vs $\varepsilon$ for the density $\rho$ $(\circ)$ (differential component) and the flux of the
momentum $z$ $(*)$ (stiff component). Top: left panel ARS(3,4,3) scheme, right panel Mod-ARS(3,4,3) scheme. Botton: BHR(5,5,3) scheme.}\label{fig2}
  \end{center}
\end{figure}

Figure \ref{fig2} represents the convergence rate of some IMEX R-K scheme
computed on a smooth test problem by grid refinement using
three different grids.  
We have obtained an improvement for the convergence of algebraic component for the Mod-ARS(4,3,4) scheme.
In fact, on the left panel we have increased the convergence rate for
sufficiently stiff parameters ($\varepsilon < 10^{-4}$). These results show a
third-order accuracy for small and large values of $\varepsilon$ and note that for intermediate
values of the parameter $\varepsilon$ ($10^{-4} < \varepsilon < 10^{-2}$) we have a slight deterioration of the
accuracy.
As it is evident in the right panel from the figure \ref{fig2} BHR(5,5,3) shows an 
almost uniform third-order accuracy in the whole range of $\varepsilon$.

\section{IMEX Runge-Kutta schemes for hyperbolic systems with diffusive relaxation}\label{Sec2}
The purpose of this section is to give a review on effective methods for the numerical
solution of hyperbolic systems with diffusive relaxation.

As the relaxation parameter vanishes,
the characteristic speeds of the system diverge, and the system reduces to a parabolic-type
equation (typically a convection-diffusion equation). 

A simple prototype of hyperbolic system with relaxation term is given by:
  \begin{align*}
    \partial_\tau u + \partial_\xi V & =  0, \\
    \partial_\tau V + \partial_\zeta p(u) & = -\frac{1}{\varepsilon}( V - Q(u))
  \end{align*}
 where $u=u(x,\tau), V=V(x,\tau) \in \RR$,
and $\varepsilon > 0$ is the relaxation time.

When looking for long time behavior of the solution of the previous system, it is more appropriate to rescale time and the variable $V$, according to the so called diffusive scaling:
  \begin{align*}
    \tau = t/\varepsilon, \quad V =  \varepsilon v, \quad
    \xi = x, \quad q(u) = Q(u)/\varepsilon, 
  \end{align*}
thus obtaining the general diffusive relaxation system given by:
\begin{eqnarray}\label{diffsyst}
\begin{array}{rcl}
  \displaystyle \partial_t u + \partial_x v  & = &  0,\\ 
  \displaystyle \partial_t v + \frac{1}{\varepsilon^2}\partial_x p(u)  & = & 
  \displaystyle -\frac{1}{\varepsilon^2} (v-q(u))  
\end{array},
\end{eqnarray}
where $p'(u) > 0$. This system is hyperbolic with two distinct real characteristics speed $\sqrt{p'(u)}/\varepsilon$.

In the small relaxation limit, $\varepsilon \rightarrow 0$ the system relaxes towards the system
  \begin{eqnarray}\label{eqlim}
  \begin{array}{rcl}
    \partial_t u + \partial_x q(u) & = & \partial_{xx} p(u),\\
    v & =& q(u) - \partial_{x} p(u).
    \end{array}
  \end{eqnarray}

The sub characteristic conditions, \cite{Liu1}, is automatically satisfied for small
$\varepsilon$ ($|q'(u)|^2 < {p'(u)}/{\varepsilon^2}$), i.e. the main stability
condition for the diffusive relaxation system. The simplest form of (\ref{diffsyst}) is to assume
$p(u) = u$ and $q(u) = 0$, then from (\ref{eqlim}) we obtain the classical heat equation $u_t = u_{xx}$. 

The attention is devoted to the construction of methods for the numerical
solution of system (\ref{diffsyst}) that are able to capture the asymptotic behavior as $\varepsilon \to 0$.
Solving (\ref{diffsyst}) numerically is challenging due to the stiffness of the problem both in the convection and in the relaxation terms.

In general, Implicit-Explicit (IMEX)
Runge-Kutta schemes represent a powerful tool for the time discretization
of stiff systems. Unfortunately, since the characteristic speed of the hyperbolic part is of order $1/\varepsilon$, standard IMEX Runge-Kutta schemes developed for hyperbolic systems with stiff relaxation \cite{ARS, CK, typeA, BR2009} fail in such parabolic scaling, because the CFL
condition would require $\Delta t = \mathcal{O}(\varepsilon \Delta x)$. 
Of course, in the diffusive regime where $\varepsilon < \Delta x$, this is very restrictive since for an explicit method a parabolic condition $\Delta t = \mathcal{O}(\Delta x^2)$ should suffice.

Most previous work on asymptotic preserving schemes for hyperbolic systems and kinetic equations
with diffusive relaxation focus on schemes which in the limit of in the finite stiffness become consistent explicit
schemes for the diffusive limit equation \cite{NP, JPT, Klar, ML}.
In those paper the authors separate the hyperbolic part into a non stiff and a stiff part and bring
the stiff part to the r.h.s., treating it implicitly. As we shall see, this can be explicitly done in several diffusive relaxation models. In all above approaches the resulting schemes, the limit scheme as $\varepsilon \to 0$ are an explicit scheme for the diffusion-like equation, with the usual parabolic CFL restriction on the time step: $\Delta t \approx \Delta x^2$.
Schemes that avoid such time step restriction and provide fully implicit solvers have been analyzed in \cite{BR1, BPR}, where a new formulation of
the problem (\ref{diffsyst}) was introduced. In the next section we review two different approaches in order to treat problem (\ref{diffsyst}) 
and some generalizations.

\subsection{Removing parabolic stiffness}
The schemes constructed with the approach outlined above converge to an explicit
scheme for the limit diffusion equation, i.e. heat equation, and therefore they are subject to the classical parabolic
CFL restriction $\Delta t \leq C\Delta x^2$. In order to overcome such a restriction we adopt a penalization
technique, based on adding two opposite terms to the first equation in (\ref{diffsyst}), and treating one explicitly
and one implicitly. 

Let us consider the simplest example of hyperbolic system
with parabolic relaxation, obtained by setting $q(u)=0$ and $p(u) = u$
in Eqs.(\ref{diffsyst}). 
By adding and subtracting the same term on the right hand side we
obtain:
\begin{eqnarray}\label{IMEX-I}
\begin{array}{l}
   u_t   = - (v +  \mu u_{x})_x + \mu u_{xx},\\
   v_t   =  -  \frac{1}{\varepsilon^2} (u_x + v).
\end{array}
\end{eqnarray}
In the first equation the term $-(v+\mu u_x)_x$ will be treated
explicitly, while the second term is treated implicitly. 
IMEX schemes based on this approach will be called 
IMEX-I, to remind that the term containing $u_x$ in the
second equation is implicit, in the sense that it appears at the new
time level. 

Notice that the term $v+u_x$ appearing in the second equation 
is formally treated implicitly, but in
practice $u_x$ is computed at the new time from the first equation, so
it can in fact explicitly computed.

The function 
$\mu: \RR^+\to [0,1]$ must be such that $\mu(0) = 1$, so that in the
limit $\varepsilon\to 0$ the quantity $(v +  \mu u_{x})_x$ vanishes. 
For $\varepsilon\ll 1$ such a quantity is very small, and so this term can be
treated explicitly. 
As $\varepsilon \to 0$ the method becomes an \emph{implicit
  scheme} for the limit equation, therefore the parabolic restriction
on the time step is removed. 

Linear stability analysis can be performed  on this simple problem, for the first order IMEX scheme, i.e.\ a backward-forward Euler method, 
both in the space continuous case, and using classical central 
differencing to approximate the first derivatives. 
For small values of $\varepsilon$ and for $\mu=1$ one obtains 
the following stability conditions (in the continuous case in space)
\[
   \begin{array}{c} 
    \displaystyle 
  {\xi^2\Delta t \leq
   \frac{1-4\varepsilon^2\xi^2}{4\varepsilon^2\xi^2}}\\[0.5cm]
   \end{array},
\]
the latter showing that there is almost no restriction for small values of
$\varepsilon$, even if we use central differences coupled with
  forward Euler time discretization.

High order extensions of this approach are possible, by using high
order IMEX (for details see \cite{BPR}). However, if we want the scheme to be accurate also in the
cases in which $\varepsilon$ is not too small, then we need to add two main
ingredients:
\begin{itemize}
\item no term should be added when not needed, i.e.\ for large enough
  values of $\varepsilon$, because in such cases the additional terms
  degrade the accuracy; this
  is obtained by letting $\mu(\varepsilon)$ decrease as $\varepsilon$
  increases. A possible choice, which is the one we use in our tests, is given by 
\[
  \mu = \exp(-\varepsilon^2/\Delta x)
\]
\item when the stabilizing effect of the dissipation vanishes, i.e. as
  $\mu \to 0$, then central differencing is no longer suitable, and one
  should adopt some upwinding; this can be obtained for example by
  blending central differencing and upwind differencing as 
\[
   D_x = (1-\mu)D_x^{\rm upw} + \mu D_x^{\rm cen}.
\]
\end{itemize}
In the IMEX-I approach, applying MOL, the diffusive system (\ref{IMEX-I}) can be written as
a ODEs system of the form
\begin{eqnarray*}
u' &=& f_1(u,v) + f_2(u),\\
{\varepsilon^2}v' &=& \ \ \ g(u,v).  
\end{eqnarray*}
where  $f_1(u,v) $ represents the discretization of the term $-\partial_x (v + \mu\partial_{x} p(u))$, 
$f_2(u)$ represents the discretization of $\mu\partial_{xx} p(u)$ and 
 $g(u,v)$ the discretization of the term $(-\partial_x p(u) - v + q(u))$. 

When $\varepsilon \to 0$ 
the solution is projected onto the manifold
$\mathcal{M} = \{(u,v) \in \mathbb{R} | g(u,v) = 0\}$.
If we assume that the equation $g(u,v) = 0$ can always be solved for $v$, and denote $v= G(u)$ the solution, then the differential variable $u$ satisfies 
\[
u' = \hat{f}_1(u) + f_2(u),
\]
with $\hat{f}_1(u) = f(u,G(u))$. The previous system is called the \emph{reduced system}. 

It would be desirable that the IMEX scheme projects the numerical solution onto the manifold 
$\mathcal{M}$ as $\veps\to 0$. In paper \cite{BPR} we proved that a sufficient condition for an IMEX scheme to project the solution onto the manifold $\mathcal{M}$ is that the scheme is {\em globally stiffly accurate}.

An implicit RK scheme is said {\em stiffly accurate} if the last row of the matrix $A$ is equal to the weights $b^T$. This ensures that the last stage is equal to the numerical solution. 
This guarantees nice stability properties of the scheme for very stiff equations (for example it ensures that the absolute stability function vanishes at infinity). 

In \cite{BPR,BR1} we extended the definition of stiff accuracy to IMEX schemes, and say that an IMEX scheme is {\em globally stiffly accurate} if the last row of both explicit and implicit RK schemes that define the IMEX are equal to the corresponding weights, i.e.\
$e_s^T A = b^T$, $e_s^T \tilde{A} = \tilde{b}^T$ with  $e_s^T = (0,...,0,1)$. 

Usually the numerical solution $(u_n, v_n)$ for all $n$ when $\varepsilon \to 0$ will not lie on the manifold $g(u, v) = 0$ since the quantity $g(u_{n}, v_{n})$ is not necessarily zero. The IMEX-I approach with a globally stuffy accurate scheme guarantees that in the limit $\varepsilon \to 0$ we obtain a globally stiffly accurate
implicit scheme and therefore  $g(u_{n}, v_{n}) = 0$ for all $n$. 

Finally in \cite{BPR} we derived additional order conditions, called algebraic conditions, that
guarantee the correct behavior of the numerical solution in the limit $\varepsilon$ maintaining the classical accuracy
in time of the scheme. We obtained such
algebraic order conditions using the classical technique by comparing
the Taylor expansion in time of the numerical solution with the one of the exact solution.
More details about this approach, as well as some rigorous analysis
can be found in \cite{BPR}.

\subsection{Additive Approach}
In the previous approach there may be a subtle difficulty when it
comes to applications, namely it is not clear how to identify the
hyperbolic part of the system, i.e.\ what is the term that should be
included in the numerical flux if I want to use my favorite shock capturing FV or
FD scheme?
We proposed in \cite{BR1} an alternative approach, in which we treat the whole hyperbolic part explicitly.
For practical applications, it would be very nice to treat the whole
term containing the flux explicitly, while reserving the implicit
treatment only to the source, according the scheme:
\[
\begin{array}{r}
u_t  \\
v_t  \\
{~}
\end{array} 
\begin{array}{c}
 = \\
 = \\
{~}
\end{array}
\begin{array}{c}
  \framebox{$
    \begin{array}{c} - v_x\\ -u_x/\varepsilon^2 
    \end{array}$} \\
  {\rm [Explicit]}
\end{array}
\hspace{2mm}
\begin{array}{c}
  \framebox{$
    \begin{array}{cl}
      &   \\
      - & v/\varepsilon^2 
    \end{array}
    $} \\
  {\rm [Implicit]}
\end{array}
\begin{array}{c}
\\[-2mm]
 \quad {\rm (Additive)}\\[2mm]
\\
\end{array}
\]

We call such an approach {\em additive} and the corresponding schemes
are denoted IMEX-E, to emphasize that the hyperbolic part is treated explicitly.

Such schemes should be easier to apply, because the fluxes retain
their original interpretation. However, the approach seems hopeless,
because of the diverging speeds.

Similarly as for the IMEX-I approach, we proposed for this approach, in order to overcome the parabolic restriction $\Delta t \approx \Delta x^2$  the same penalization technique based on adding two opposite terms to the first equation, and treating one
explicitly and one implicitly.
 
In this paper, the authors concentrated on developing IMEX R-K schemes of
type A, since they are easier to analyze with respect to the other types. They started
the analysis by introducing a property which is important in order to guarantee the
asymptotic preserving property, i.e. the scheme possesses the
correct zero-relaxation limit, in the sense that the numerical scheme applied to system (\ref{diffsyst}) should be a
consistent and stable scheme for the limit system (\ref{eqlim}) as the parameter $\varepsilon$ approaches zero independently
of the discretization parameters. IMEX R-K schemes that satisfy this property are globally stiffly accurate schemes.
Several results and a rigorous analysis about that can be found in \cite{BR1}. Most numerical
tests are reported in \cite{BR1} for IMEX-E approach
and the results are compared with those obtained by other methods available in the
literature.

\section{Applications.} This section is devoted to the presentation of some applications of the previous two approaches for the treatment of hyperbolic systems with diffusive relaxation.

\subsection{Kawashima-LeFloch's nonlinear relaxation model} 
\begin{figure}
\begin{center}
\includegraphics[width=0.3\textwidth]{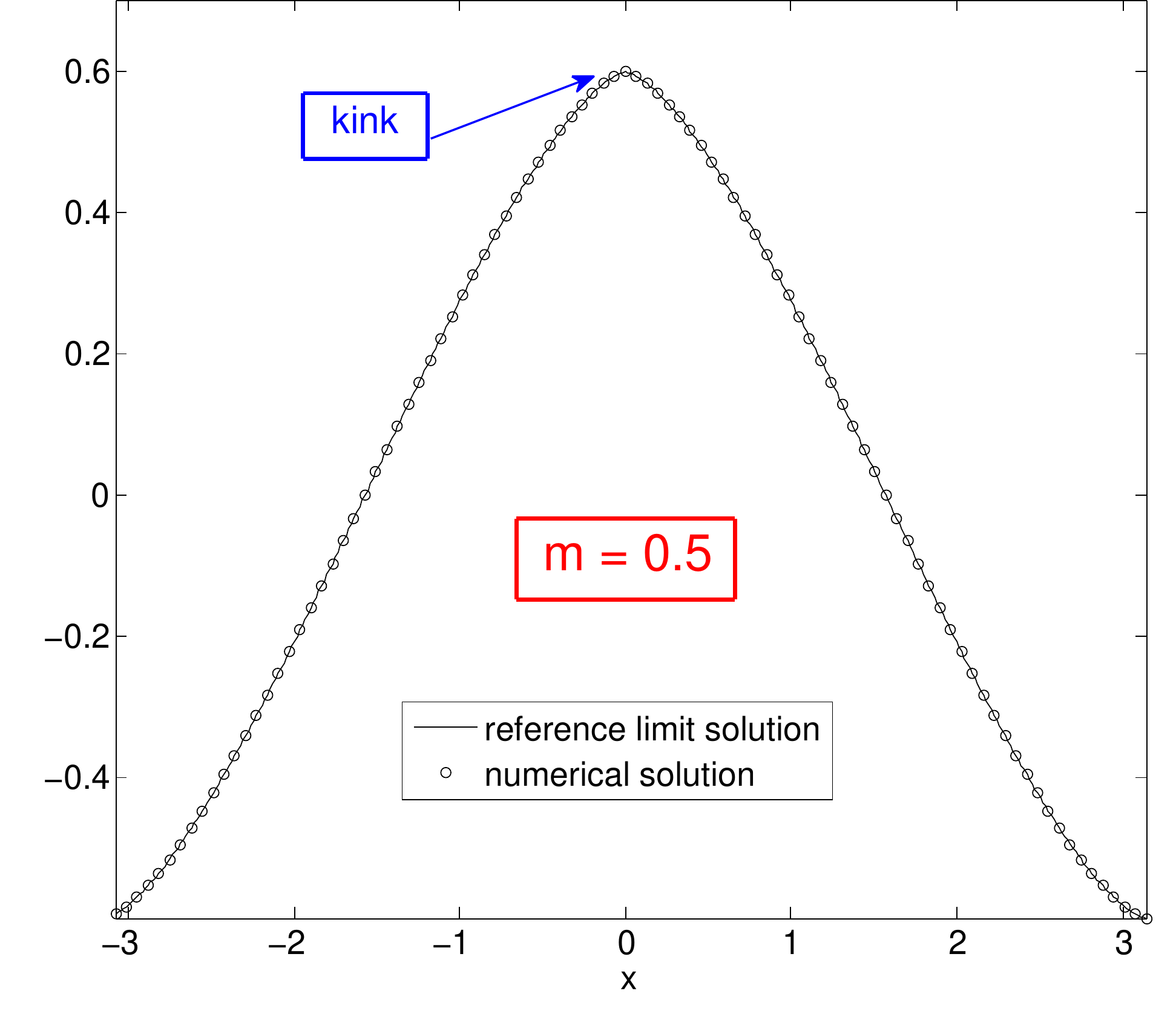}
\includegraphics[width=0.3\textwidth,height=0.265\textwidth]{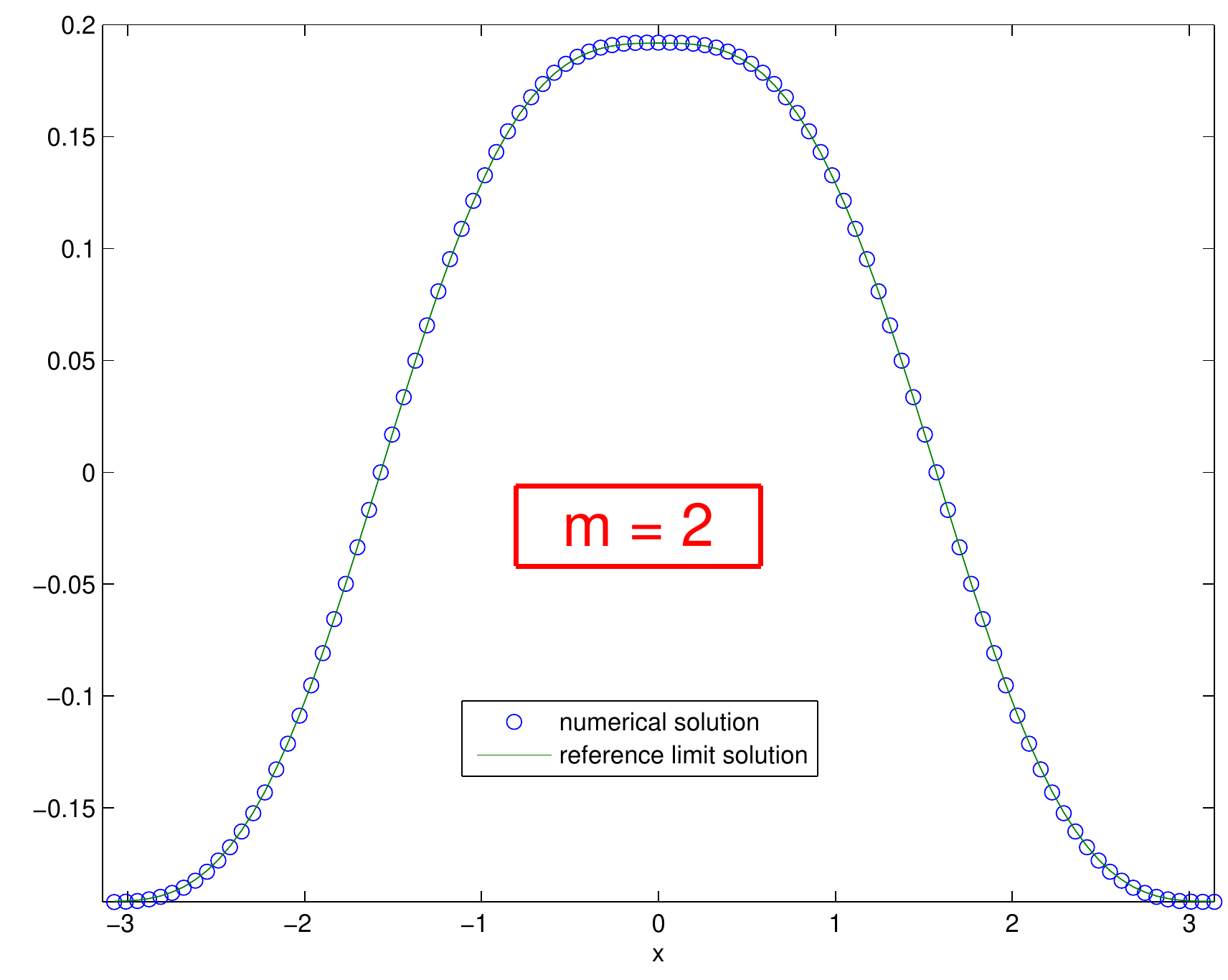}
\caption{Numerical solution with $N = 96$ cells. Solid line: reference
limit solution with $N = 384$ cells at time $T = 1$.
IMEX-E approach, $\veps = 10^{-4}$ and $\Delta t = C \Delta x^2$. $u(x,0) = \cos(x)$, $v(x,0) = \sin(x)$. 
Left $C = 1$.  Right $C = 0.025$.}\label{regular}
\end{center}
\end{figure}

\begin{figure}
\begin{center}
\includegraphics[width=0.7\textwidth]{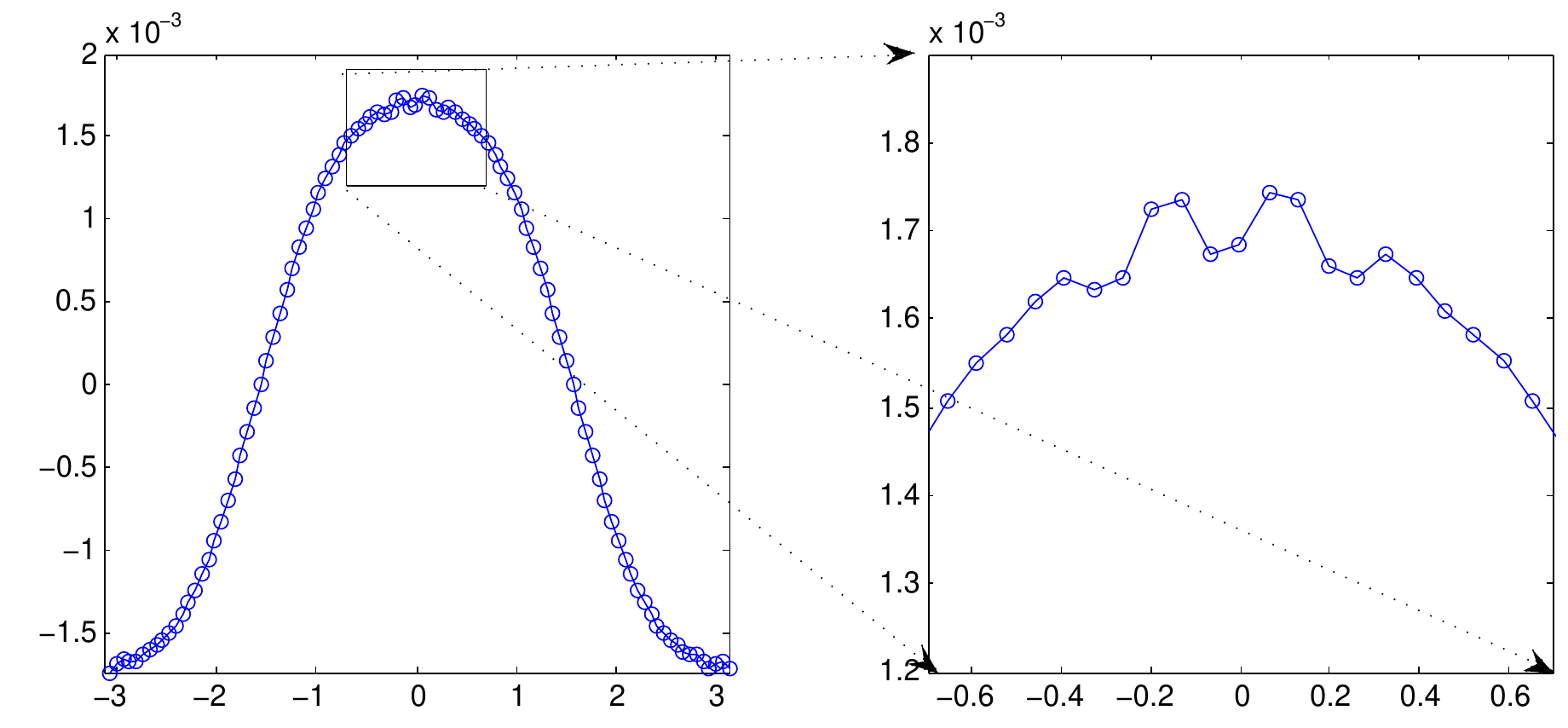}
\caption{Instabilities for $m=2$ ($\alpha = -1$), $\Delta t = C \Delta x^2$, $\veps = 10^{-4}$ $N = 96$.}\label{osc}
\end{center}
\end{figure}

Fully nonlinear relaxation terms arise, for instance in presence of nonlinear
friction and, in this section we want numerically study the following non-linear relaxation model, first introduced by Kawashima and LeFloch \cite{KLF}, i.e.
\begin{eqnarray}\label{KLF}
\begin{array}{c}
u_t + v_x  = 0,\\
\veps^2 \, v_t + b(u)_x = - |v|^{m-1} \, v + q(u).
\end{array}
\end{eqnarray}
Provided $b'(u) > 0$, system (\ref{KLF}) is strictly hyperbolic system of balance laws.
In the stiff relaxation $\varepsilon$, ($\varepsilon \to 0$) we have
\begin{align*}
u_t  &= (|-q(u) + b(u)_x|^{\alpha}(-q(u) + b(u)_x) )_x,\\
|v|^{m-1} \, v & =  q(u) - b(u)_x,
\end{align*}
which is a fully nonlinear parabolic equation in $u$ with $\alpha = -1 + 1/m$. We distinguish between
\begin{align*}
\rm{sub-linear} & : 0 < m < 1,\\
\rm{linear} & :   m =1,\\
\rm{super-linear} & : m > 1.
\end{align*}
In its simplest form we assume $b(u) = u$, $q(u) = 0$ and we get:
\begin{align*}
u_t + v_x & = 0,\\
\veps^2 \, v_t + u_x & = - |v|^{m-1} \, v.
\end{align*}
As $\veps\to 0$ this relaxes to 
\begin{eqnarray}\label{KF2}
u_t  =  \left( |u_x |^{\alpha} u_x  \right)_x,\qquad |v|^{m-1} \, v = - u_x.
\end{eqnarray} 
Very interesting cases are both $m<1 \> (\alpha>0)$ and $m\ge1 \> (\alpha \le 0)$,
The profile of the solution computed with $N = 96$ points is reported in Fig. \ref{regular}. 
But by integrating for a longer time, the nonlinear parabolic equation (\ref{KLF}) has regular solutions if $m > 1$, i.e. $\alpha \le 0$, while it develops singularities in the derivatives if $0 < m < 1$, i.e.  $\alpha > 0$.
In fact, for $m=2$ ($\alpha = -1$) integrating for a longer time  $T=1.77$, some instabilities appear (see Figure \ref{osc}).
The reason of such instabilities is that equation (\ref{KF2}): 
\[
u_t  =  \left((\alpha + 1)|u_x |^{\alpha}\right)u_{xx},
\]
where the non-linear diffusion coefficient $\nu$ is
\[
   \nu = (1+\alpha)|u_x|^\alpha
\]
which suggests the following condition in the nonlinear case
\begin{eqnarray}\label{stab}
  (1+\alpha)|u_x|^\alpha \frac{\Delta t}{\Delta x^2} \le 1
\end{eqnarray}
but the equation  (\ref{KF2}) diverges near local extrema when $\alpha<0$ ($m>1$).
This condition (\ref{stab}) is used to determine the optimal time step for $m \le 1$, no time step can guarantee stability near local extrema if $m > 1$. 
In \cite{BFR}, the same penalization technique proposed in \cite{BPR, BR1} in order to remove the parabolic stability restriction has been used. 

Then we write the system in the form
\begin{align*}
 u_t          & = - (v+\mu(\veps) |u_x|^{\alpha}u_x)_x + \mu(\veps)(|u_x|^{\alpha}u_{x})_x   \\
 \veps^2 v_t   & =  -u_x-|v|^{m-1}v.
\end{align*}
Now in order to treat this system by the IMEX-I or IMEX-E approach this requires that the term $(|u_x|^{\alpha}u_{x})_x$ is treated
implicitly. But some difficulty arises, in fact, when $\varepsilon \to 0$, the limit equation is non-linear parabolic and 
fully implicit would be very expensive. 

\begin{figure}
    \begin{center}
    \includegraphics[width=0.5\textwidth]
    {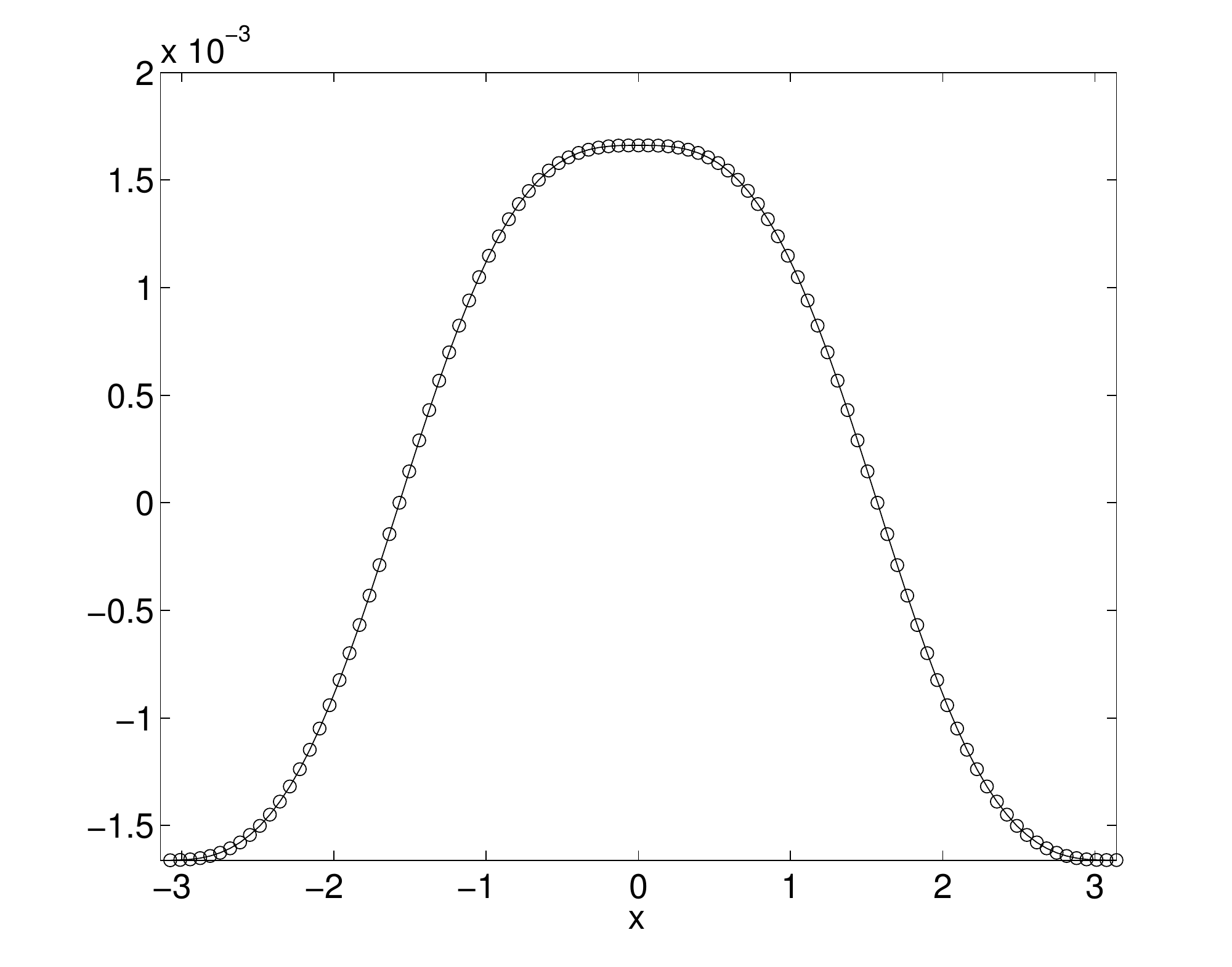}
    \caption{Numerical solution with $N = 96$ cells at time $T = 1.77$ for $m = 2$, $\varepsilon = 10^{-4}$ and
    $\Delta t = C \Delta x$, with $C = 0.25$.}\label{LongT}
   \end{center}
\end{figure}

In \cite{BFR} a new approach has been used in order to solve the term $(|u_x|^{\alpha}u_{x})_x$ where a very efficient method for the numerical solution of such an equation has been introduced. Indeed the idea is to write the equation as a system as
\begin{eqnarray}
     y' = F(y^{\ast},y)  
     \end{eqnarray}
with $F$ function non-stiff in the first variable and stiff in the second one. 
To be more specific,  in our case $F(y^{\ast}, y)$ is given by 
$
y = \left(\begin{array}{l}
u\\
v
\end{array}\right),\,
y^{\ast} = \left(\begin{array}{l}
u^{\ast}\\
v^{\ast}
\end{array}\right),
$ and 
$$
F(y^{\ast}, y) = \left(\begin{array}{l}
 -(v^{\ast}_x + \mu(\varepsilon)(|u^{\ast}_x|^{\alpha} u^{\ast}_x)_x) +  \mu(\varepsilon)(|u^{\ast}_x|^{\alpha}u_x)_x \\
-u_x + |v|^{m-1}v 
\end{array}\right).
$$
Additive RK for this class of problems can be constructed, in particular we showed that in order to compute the numerical solution we need to require that $b_i=\tilde{b}_i$ for $i$ (see \cite{BFR}), then a good choice is to consider IMEX-I approach, whereas IMEX-E approach requires that the Runge-Kutta IMEX is globally stiffly accurate, i.e. $\tilde{b}_i \neq b_i$ for all $i$ \cite{BR1}. 
Using this new approach one can solve the relaxation system without
parabolic CFL, i.e. $\Delta t = 0.25 \Delta x$ and $T= 1.77$, Fig. \ref{LongT}.
   \begin{itemize}
   \item Time step is about 150 times larger than in the explicit method.
    \item The case $m = 2$, i.e. $\alpha = -1$, we set a \emph{TOL}
    for computing $(|u|_x + TOL)^{\alpha}$, in order to avoid that the derivatives goes to infinity.
\end{itemize}

\subsection{R13: a regularized Grad's 13 moment method} 
\label{R13}
Grad's moment method is a technique used to close the infinite hierarchy of moments arising from the Boltzmann equation or rarefied gases.
It is an example of hyperbolic relaxation: the Boltzmann equation relaxes to the hyperbolic system of Grad's equations.
Sometimes parabolic systems provide more accurate physical description (e.g.\ Navier-Stokes equations are very successful in practice, although they are not hyperbolic). Some researchers, mainly Manuel Torrilhon and Henning Struchtrup \cite{TS} developed a parabolic extension of Grad's approach, called R13. When derived from the Boltzmann equation, this can be viewed as a parabolic relaxation.

In this section we present some results for the asymptotic accuracy for boundary value problems,
which emerges from a 1D simplification of the R13 system that describes a Poiseuille-flow \cite{To}. The system takes the form
\begin{equation}
U_\tau + F(U)_\xi =  - \frac{1}{\veps}P(U) + G \label{eqPF}
\end{equation}
Here, the variables are $U = (u, v, w)^T$ with velocity $u$, shear stress $v$ and parallel heat $w$.
Furthermore, we have
\begin{eqnarray}\label{AA}
F(U) =  AU, \quad 
A = 
\left(\begin{array}{ccc}
0 & 1 & 0\\
1/2 & 0& 1/2\\
0 & 1 & 0
\end{array}\right), \quad
P(U) = \left( \begin{array}{c}
0\\
v\\
w
\end{array}\right), \quad
G = \left( \begin{array}{c}
g\\
0\\
0
\end{array}
\right)
\end{eqnarray}
where here the parameters $g$ and $\varepsilon$ are the external force and the relaxation time.
Explicitly, we write system (\ref{eqPF}) as
\begin{eqnarray}\label{expPF}
u_\tau + v_\xi &=& g, \nonumber        \\
v_\tau +\frac12{(u + w)_\xi} &=& -\frac{v}{\varepsilon},\\
w_\tau +  v_\xi &=& -\frac{w}{\varepsilon}.\nonumber
\end{eqnarray}
We consider a bounded domain $\xi \in [-1, 1]$ where we have to prescribe
boundary conditions.
In \cite{To}, the authors used the following boundary conditions for $v$: 
\begin{eqnarray}\label{BC1}
\begin{array}{c}
v|_{\xi = \pm 1} = \pm(\alpha u + \beta w)_{\xi = \pm 1},         
\end{array}
\end{eqnarray}
with $\alpha> \beta > 0$  some parameters.  
In the numerical experiments we chose the following values, $g = 1$,
$\alpha = 0.7$, $\beta = 0.3$. A steady state solution for system
(\ref{eqPF}) is given by
\begin{eqnarray}\label{stst}
u_s(\xi) = g\left(\frac{1+ \varepsilon \beta}{\alpha} + \frac{1}{\varepsilon}(1-\xi^2)\right), \ \ v_s(\xi) = g\xi, \ \ w_s(\xi) = -\varepsilon g.
\end{eqnarray}
We consider numerical tests whose solution converges to such steady state.
We note that as we use high order reconstruction for the fuxes, then we need two layers of ghost cells that can be
obtained using the boundary values. This part of the discretization is most important, because
the efficiency  of the whole method heavily depends on the choice of the correct boundary values
and extrapolation methods.

Now we will focus our attention to the following system
\begin{eqnarray}\label{diff}
\tilde{u}_t + v_x &=& g,\nonumber        \\
v_t +\frac{1}{2}(\frac{\tilde{u}}{\varepsilon^2} + \tilde{w})_x &=& -\frac{v}{\varepsilon^2},\\
w_t +  \frac{v_x}{\varepsilon^2} &=& -\frac{\tilde{w}}{\varepsilon^2}. \nonumber        
\end{eqnarray}
obtained by (\ref{expPF}) under the diffusive scaling $t=\varepsilon \tau$, $x = \xi$, $\tilde{u} = \varepsilon u$ and $\tilde{w} = w/\varepsilon$.

Concerning the space discretization we consider a finite volume discretization as done in \cite{To}. In our diffusive approach the matrix $A$ in (\ref{AA}) has the following expression
\begin{equation}
   A = 
   \left(\begin{array}{ccc}
       0 & 1 & 0\\
       1/2\varepsilon^2 & 0 & 1/2\\
       0 & 1/\varepsilon^2 & 0\\
   \end{array}\right) 
\end{equation}

In the small relaxation (or diffusion) limit, i.e. when $\varepsilon \rightarrow 0$, the behaviour of the solution to (\ref{diff}) is governed by 
\begin{eqnarray}
\tilde{w} = - v_x,\ \ \ \ v = \frac{-\tilde{u}_x}{2},
\label{wv}
\end{eqnarray}
and
\begin{eqnarray}\label{eq1}
\tilde{u}_t  = \frac{\tilde{u}_{xx}}{2} + g.
\end{eqnarray}  
Now consider boundary conditions for (\ref{diff}) which are consistent
to the limit system (\ref{wv},\ref{eq1}).
%

\subsubsection{Boundary Treatment} In this section we derive boundary
conditions which are in agreement with the stationary solution. 

From the steady state condition of Eq.~(\ref{diff}) we get
\begin{eqnarray*}
v_x & = & g,\\
\frac{\tilde{u}_x/\varepsilon^2 +\tilde{w}_x}{2} & = & -\frac{v}{\varepsilon^2}.\\
v_x & = & -\tilde{w}
\end{eqnarray*}
We observe that compatibility with stationary solutions implies:
\begin{eqnarray}
\displaystyle 	\tilde{w}|_{\pm 1} & = & -g,\label{1}\\
\displaystyle	 v_x|_{\pm 1} & = & g,\label{2}\\
\displaystyle	 \left(\tilde{u}_x + \varepsilon^2 \tilde{w}_x
\right)|_{\pm 1} & = & -2v|_{\pm 1} \label{3}.
\end{eqnarray}
Such conditions are compatible with condition
\begin{eqnarray}\label{4}
\begin{array}{l}
\displaystyle \tilde{u} |_{ \pm 1} =  \pm  \frac{\epsilon v \mp \beta \tilde{w} \epsilon^2}{\alpha}.\\
\end{array}
\end{eqnarray}
for the stationary solution (\ref{stst}).
Therefore one can solve the system with boundary conditions (\ref{1}), (\ref{2}), (\ref{3}) or (\ref{1}), (\ref{2}) and (\ref{4}). 
In both cases one obtains convergence to the stationary solution. 

System (\ref{eqPF}) is discretized by second order finite volume for
the internal points. Ghost points are used out of the boundary to
impose boundary conditions. Such ghost points are computed by
extrapolation. 
For instance for the calculation of the boundary values considering (\ref{1}), (\ref{2}) and (\ref{4}),  we can write
\begin{eqnarray}
\tilde{w}^W_{1/2} & = & -g,\\
\tilde{w}_0 & = & (8w^{W}_{1/2}-6\tilde{w}_1 + \tilde{w}_2)/3,\\
v_0 & = & v_1- g \Delta x,\\
v^W_{1/2} & = & \frac{3}{8} v_0 + \frac{3}{4} v_1 - \frac{1}{8} v_2,\\
\tilde{u}^W_{1/2} & = & -\varepsilon(v^W_{1/2} + \varepsilon \beta \tilde{w}^W_{1/2})/\alpha,\\
u_0 & = & (8\tilde{u}^{W}_{1/2}-6\tilde{u}_1 + \tilde{u}_2)/3,\\
U_{-1} & = & 3 U_0 - 3U_1 + U_2,
\end{eqnarray}
where $U =(\tilde{u},v,\tilde{w}) $. We do the same for the other part of the wall $x_{N+1/2}$.

We remark that we can improve the order of the extrapolation to the ghost cells by the following considerations. We consider the Lagrange polynomial 
\begin{eqnarray*}
L_n(x;U) = \sum_{i = 0}^n U_i \ell_i(x)
\end{eqnarray*} 
where $U=(\tilde{u},v,\tilde{w})$ and 
\begin{enumerate}
	\item $\tilde{w}_{0} \ell_0(x_{1/2}) = -g-\sum_{i=1}^n \tilde{w}_i \ell_i(x_{1/2})$,\\
	\item $v_{0} \ell^{'}_0(x_{1/2}) = g-\sum_{i=1}^n v_i \ell^{'}_i(x_{1/2})$,\\
\item $\tilde{u}_{0} \ell_0(x_{1/2}) = -\frac{\varepsilon}{\alpha}(v(x_{1/2})+\beta\tilde{w}(x_{1/2})\varepsilon)-\sum_{i=1}^n \tilde{u}_i \ell_i(x_{1/2})$.
\end{enumerate}
and we can compute 
\begin{eqnarray*}
U_k = \sum_{i=0}^{n}U_i \ell(x_k), \ \ \ k = -1,-2,...
\end{eqnarray*}
Similarly for the other side of the wall.

We remark that we tested this approach performing also a numerical simulation setting different initial conditions, i.e. introducing a little perturbations to the initial data, and we observed that after a short time the numerical solution converge to the stationary solution.  
	
\subsubsection{Removing parabolic stiffness}\label{RPS}
We rewrite system (\ref{diff}) in the following form 
\begin{eqnarray}\label{diff2}
\tilde{u}_t  &=& - v_x - \underbrace{\mu(\varepsilon)\frac{\tilde{u}_{xx}}{2} + \mu(\varepsilon)\frac{\tilde{u}_{xx}}{2}} + g,\nonumber        \\
v_t  &=& - \frac{\tilde{w}_x}{2} -\frac{1}{2}\frac{\tilde{u}_x}{\varepsilon^2}  -\frac{v}{\varepsilon^2},\\
\tilde{w}_t &=& -  \frac{v_x}{\varepsilon^2} -\frac{\tilde{w}}{\varepsilon^2}. \nonumber        
\end{eqnarray}
where we added and subtracted the term $\mu(\varepsilon){\tilde{u}_{xx}}/{2}$ in order to overcome the stability restriction that usually we have for hyperbolic system with diffusive relaxation. Here $\mu(\varepsilon)$ is such that $\mu : \mathbb{R}^{+} \rightarrow [0, \ 1]$ and $\mu(0)  =1$. When $\varepsilon$ is not small there is no reason to add and subtract the term $\mu(\varepsilon) u_{xx}$, therefore $\mu(\varepsilon)$ will be small in such a regime, i.e.\ $\mu(\varepsilon) \approx 0$. For a detailed analysis on this topic we report to \cite{BPR}.
Furthermore this reformulation allows us to design a class of IMEX Runge-Kutta schemes that work with high order accuracy in time in the zero-diffusion limit, i.e.\ when $\varepsilon$ is very small, and in a wide range of the parameter $\varepsilon$ such that the scheme maintains the accuracy uniformly for each value of $\varepsilon$.
Now we want to apply an IMEX Runge-Kutta scheme with these features to
this system considering IMEX-I approach, \cite{BPR}. In our numerical
test we consider the stiffly accurate IMEX-SSP2(3,3,2) which satisfies
all the conditions described above. 

Then, by  (\ref{diff2}), we treat the quantities 
\begin{eqnarray}
(- v_x -\mu(\varepsilon)\frac{u_{xx}}{2},- \frac{\tilde{w}_x}{2}, 0)^T
\end{eqnarray}
 explicitly and 
 \begin{eqnarray}
 (\mu{(\varepsilon)}\frac{u_{xx}}{2} + g, -\frac{1}{2}\frac{\tilde{u}}{\varepsilon^2}  -\frac{v}{\varepsilon^2}, -\frac{\tilde{v}_x}{\varepsilon^2} -\frac{\tilde{w}}{\varepsilon^2})^T
 \end{eqnarray}
  implicitly, respectively.
\subsubsection{Convergence Results}
In order to ensure the second order convergence for the
IMEX-SSP(3,3,2) scheme with the previous boundary conditions proposed,
we simulate the same periodic test case proposed in \cite{To}. We
chose $g = 0$ and the initial conditions are $u = \sin(\pi x) + 0.5
\sin(5\pi x)$, $v = 0$, $w = 0$. We simulate until $t_{end} =
\varepsilon \tau_{end}$ with $\tau_{end} = 4$, $\varepsilon
=0.01$. 

Numerical convergence rate is calculated by the formula 
\begin{eqnarray}\label{CR}
p = \log_3(E_{\Delta t_1}/E_{\Delta t_2}),
\end{eqnarray}
where $E_{\Delta t_1}$ and $E_{\Delta t_2}$ are the global errors
associated to time steps $\Delta t_1$ and $\Delta t_2$, respectively. 
$E_{\Delta t_1}$ is obtained by comparing a solution with $N=50$ with
a solution obtained using $N=150$ points, while for $E_{\Delta t_2}$
we use two solutions obtained, respectively, with $N=150$ and $N=450$
points. The number of points is tripled each time, because in this way
it is easier to compare solutions in the same location using finite
volume discretization.
 In Table \ref{conT} we show that a second order is reached for IMEX-SSP(3,3,2) scheme for all three components.

We note that we have obtained these convergence results considering
the system (\ref{diff}) without adding and subtracting any term. It is
clear from the previous considerations that a time step $\Delta t =
\mathcal{O}(\Delta x^2)$ must be chosen. It is
possible to obtain similar results considering the reformulated 
system (\ref{diff2}) and choosing a time step $\Delta t =
\mathcal{O}(\Delta x)$, although a special care has to be taken when
imposing boundary conditions in the implicit step.
\begin{table}
\begin{center}
{\begin{tabular}{|cccc|}
\hline
\lower.3ex\hbox{$N$} &  \lower.3ex\hbox{$Error_u$} & \lower.3ex\hbox{$Error_v$}  & \lower.3ex\hbox{$Error_w$}  \\
\hline
\hline
\lower.3ex\hbox{ $\ 50$} & -- & -- &-- \\
\lower.3ex\hbox{ $150$}& \lower.3ex\hbox{$8.062e-04$} & \lower.3ex\hbox{$2.530e-03$}& \lower.3ex\hbox{$1.089e-02$}\\
 \lower.3ex\hbox{$\ 450$}  & \lower.3ex\hbox{$7.838e-05$} &  \lower.3ex\hbox{$2.879e-04$} & \lower.3ex\hbox{$1.162e-03$} \\
                          \hline
                          \hline
 \lower.3ex\hbox{Order}   & \lower.3ex\hbox{$2.121$} &  \lower.3ex\hbox{$1.978$} & \lower.3ex\hbox{$2.036$} \\                       
                          \hline
  \end{tabular} }
  \end{center}  
  \caption{}\label{conT}             
\end{table}

We now investigate numerically the convergence rate for a wide range
of $\varepsilon$ considering system (\ref{diff2})  and choosing a time
step $\Delta t = \mathcal{O}(\Delta x)$. To this aim we consider the
previous test problem with the second order IMEX-SSP(3,3,2) scheme introduced before. 
Numerical convergence rate is calculated by (\ref{CR}) and time step $\Delta t = 0.3 \Delta x$.
We simulate until $t_{end} = 1$.

\begin{figure}
\begin{center} 
\includegraphics[width=0.75\textwidth]{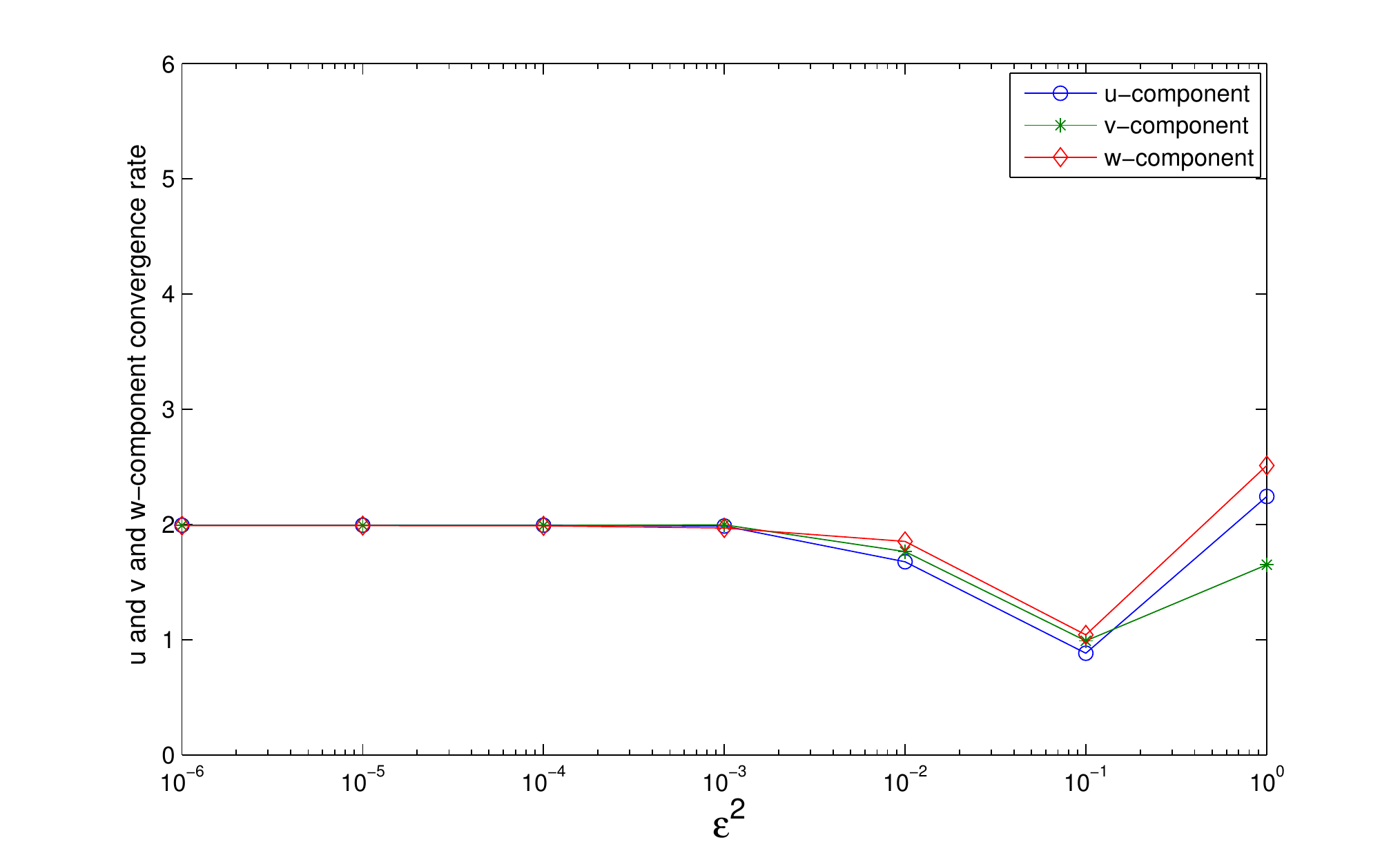}
\caption{Convergence rate for the $u$, $v$, and $w$ component versus $\epsilon^2$}\label{fig:1bis}
\end{center}
\end{figure}

Figure \ref{fig:1bis} shows the convergence rates as a function of
$\varepsilon^2$ using different values of $\varepsilon$ ranging from
$10^{-6} $ to $1$. The second order scheme tested has the prescribed
order of accuracy uniformly in $\varepsilon^2$  until $\varepsilon$ is
small. Instead, for values of $\varepsilon$  large, say $10^{-1}$, a
degradation of accuracy is observed. This phenomenon requires further investigation as mentioned in \cite{BPR}.

\subsubsection{Convergence to the steady state solutions}\label{NumCon}
In this numerical test we show how starting form arbitrary initial
conditions and considering the stiffly accurate IMEX-SSP2(3,3,2), the
IMEX-I approach proposed in section (\ref{RPS}), (see for details
\cite{BPR}), provides a numerical solution that converges to the
steady state solution (\ref{stst}) in a number of time steps much
smaller than the one needed by classical IMEX methods. 

We consider $g = 1$, $\alpha = 0.7$, $\beta = 0.3$ and we choose
$\varepsilon = 10^{-4}$ (\emph{diffusive regime}). The final time is
$\tau = 10$, the domain is $I = \{x: x \in [-1, 1]\}$ and 
$\Delta t_H = 2.5 \Delta x$ with $N = 50$ grid points. This $CFL$ number has been empirically
adjusted to approximate the largest one that maintains stability.  As initial data we consider
\begin{equation}
u_0 = \frac{\varepsilon}{\alpha}\left( (C + \beta \varepsilon)x  - g\right)\quad
v_0 = gx + C\quad
w_0 = -x^2.
\end{equation}
This initial conditions are compatible with the boundary conditions
(\ref{1}), (\ref{2}), (\ref{4}). We plot the numerical solution at
different final times $0.5, \ 1, \ 1.5, \ 3$ and $10$. At the final
time the numerical
solution is in perfect agreement with the steady state solution
(\ref{stst}) after $100$ time steps.
We remark that the steady state solution
is in practice reached a smaller time, say $t=5$. 
We chose a long time in order
to show that the numerical solution reaches the steady state with no 
oscillations.  IMEX-I approach with the penalization technique
described in 
Sec.~\ref{RPS} allows a time step $\Delta t$ with a hyperbolic
stability restriction rather than the parabolic one
typical of explicit
schemes for diffusion problems. 
Indeed, if we compute the numerical solutions $u$, $v$ and $w$
of system (\ref{diff}) without adopting the penalization technique,
when $\varepsilon$ is very small a stability parabolic
restriction like $\Delta t_P = CFL \Delta x^2$ is required because the
IMEX R-K method becomes an explicit one in the limit case $\varepsilon
\to 0$. In this case we consider $CFL = 2.5$ and we note that thanks
to the better stability properties of the new approach, the time
step $\Delta t_H$ is about 25 times bigger then $\Delta t_P$.
\begin{figure}
\centering 
\includegraphics[width=0.6\textwidth]{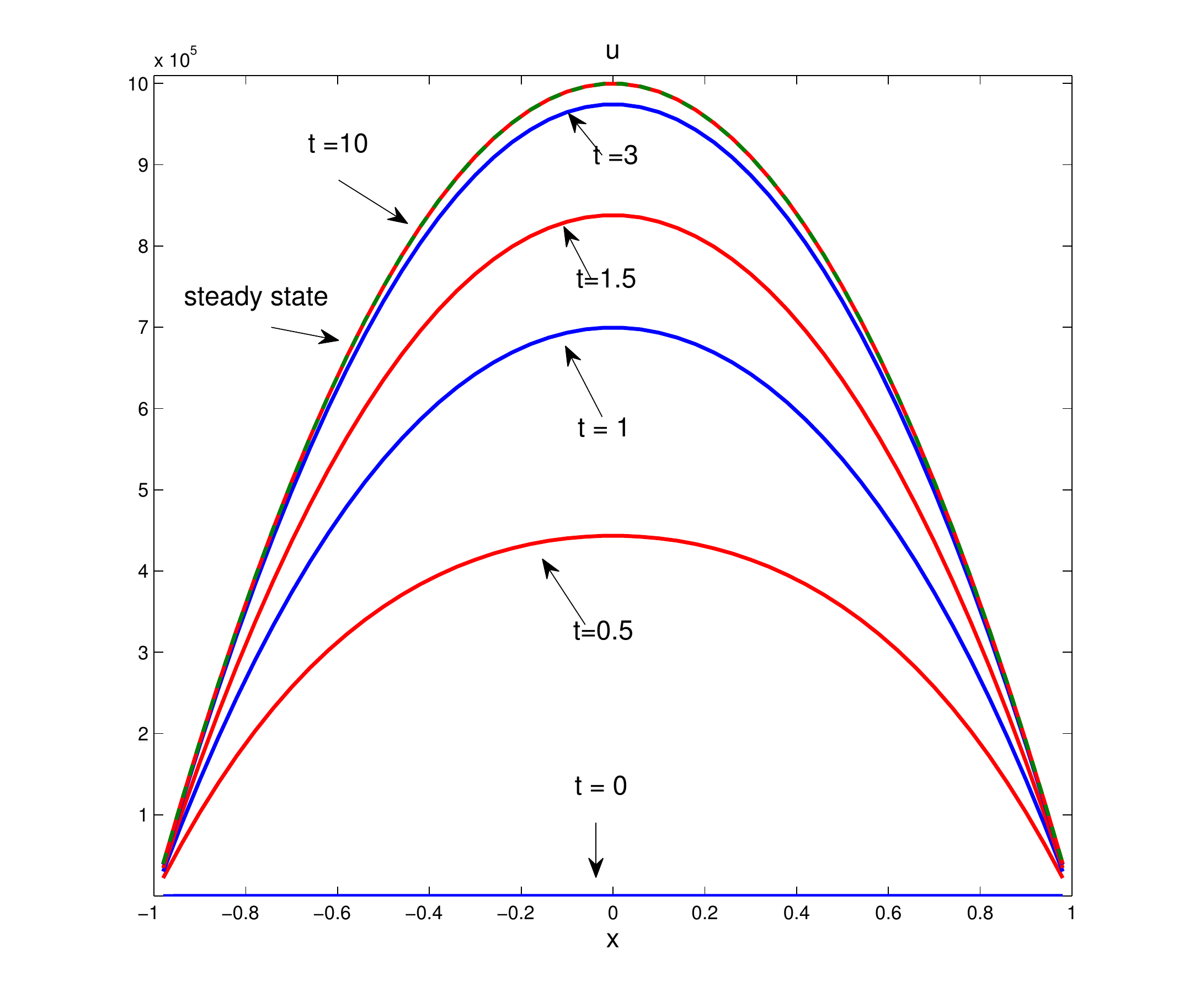}
%
\includegraphics[width=0.6\textwidth]{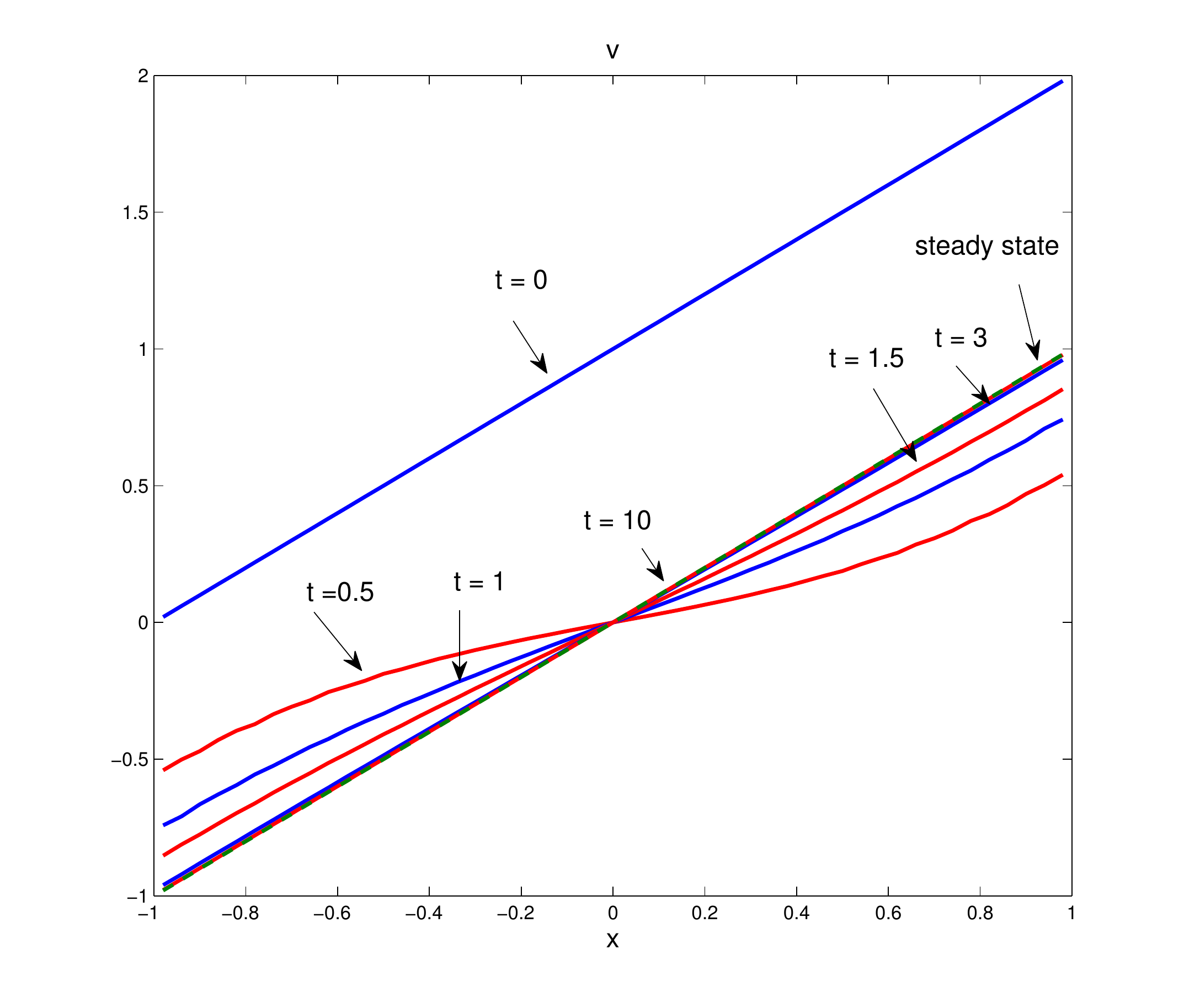}
%
\includegraphics[width=0.6\textwidth]{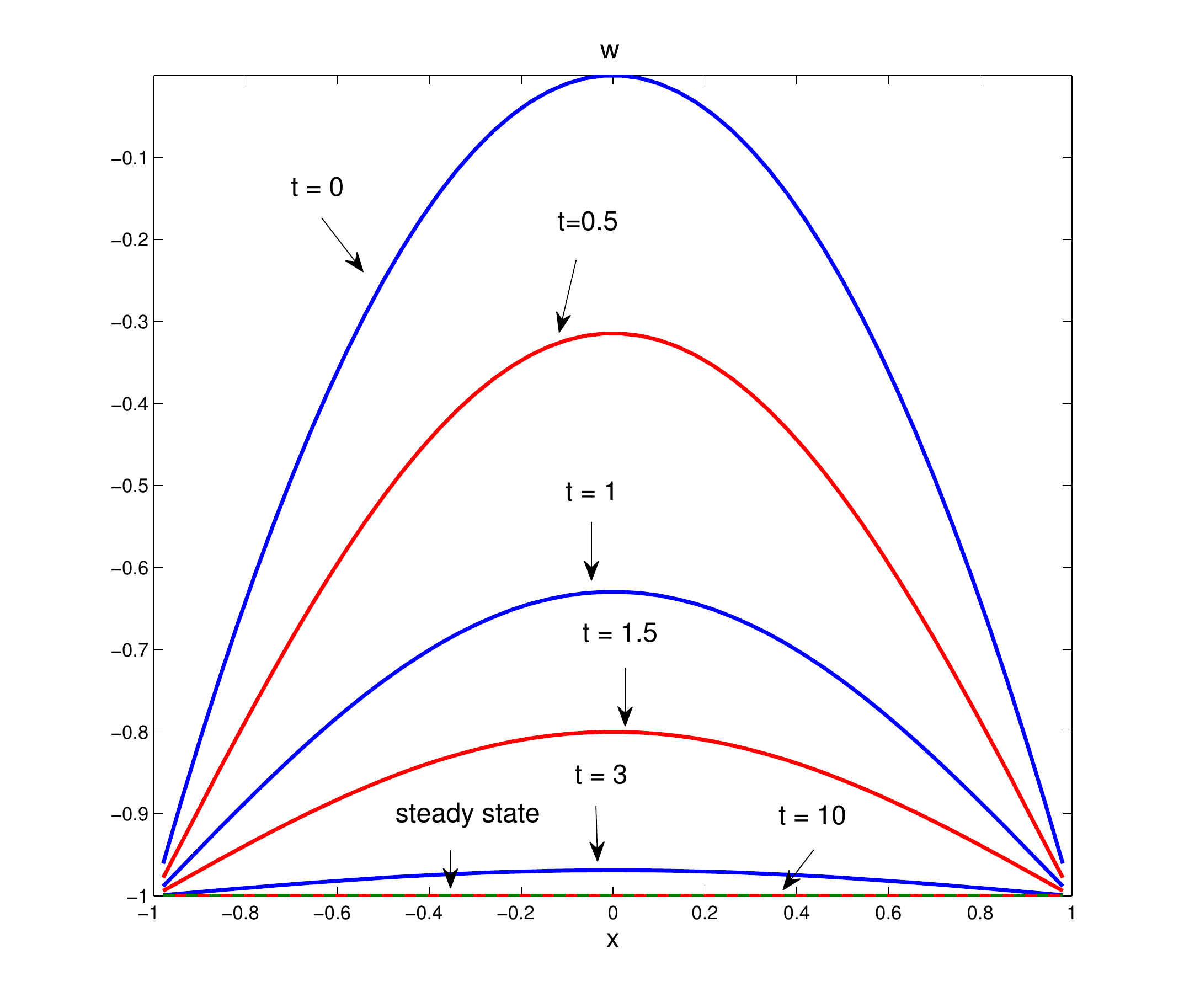}
\caption{Convergence to the steady state for the R13 model
  problem. From top to bottom: $u$, $v$, and $w$ profiles at different
  times. Number of grid points $N=50$. Time step $\Delta t_H =
  2.5\Delta x$.
}
\end{figure}

\section{Conclusions}
We gave a brief review of modern IMEX Runge-Kutta schemes for
hyperbolic systems in presence of stiff relaxation. Both hyperbolic
and parabolic relaxations are considered, in the framework of
conservative finite difference space discretization, which is the
simplest approach to construct high order shock capturing schemes for
such problems. 

In the hyperbolic relaxation case, most IMEX schemes in the literature
are able to capture the correct relaxed limit, converging to explicit
schemes for the relaxed system. If high accuracy is required for a
wide range of values of the relaxation parameter, then suitable
conditions have to be imposed on the coefficients of the scheme in
order to guarantee uniform accuracy, based on the analysis developed
in \cite{Boscarino2007}. 

The parabolic case is more subtle, since the characteristic speeds of
the hyperbolic part diverge as the stiffness parameter vanishes. 
Numerical schemes commonly found in the literature for this family of
problems converge to an explicit scheme for the limit parabolic
equation, thus requiring a parabolic type CFL restriction on the time
step. Recently developed schemes overcome such problem, using a
penalization technique, and providing IMEX
schemes that relax to an implicit scheme for the limit diffusion
equation, of to an IMEX scheme for the limit convection-diffusion
equation (according to the form of the relaxation term)
\cite{BPR,BR1}. 
Suitable modification of such schemes can be adapted to problems that
relax to genuinely non-linear diffusion equations \cite{BFR}.
IMEX schemes with the penalization techniques are applied here to a
model problem coming from Extended Thermodynamics, providing a much
more efficient tool to solve the problem with a number of time steps
considerably smaller than the one required by other schemes present in
the literature. 

Several open problems remain. In particular we mention two problems
that may attract the attention of researchers in this area. The first
one is the extension of the uniform accuracy analysis performed in the
case of hyperbolic relaxation to the more difficult problem of the
parabolic relaxation. The second problem consists in exploiting the
stabilization effect of the penalization technique adopted to improve
the stability properties of the IMEX schemes for the parabolic
relaxation to more a more general framework, extending the work
already performed in \cite{smereka} and \cite{jin-filbet} in specific
cases.

\end{document}